\documentclass[11pt,reqno]{amsart}
\usepackage{amssymb,amscd,amsbsy,mathrsfs}

\usepackage[cp1251]{inputenc}

\newcommand{\lb}{\linebreak}
\newcommand\Li{{\rm Lip}}
\newcommand{\bS}{{\boldsymbol S}}
\newcommand{\bs}{\boldsymbol}
\newcommand{\bay}{\begin{eqnarray}}
\newcommand{\ey}{\end{eqnarray}}
\newcommand{\trace}{\operatorname{trace}}
\newcommand{\R}{{\Bbb R}}
\newcommand{\const}{\operatorname{const}}
\newcommand{\rf}[1]{(\ref{#1})}
\newcommand{\be}{\infty}
\newcommand{\m}{{\boldsymbol m}}
\newcommand{\e}{\varepsilon}
\renewcommand{\o}{\omega}
\newcommand{\supp}{\operatorname{supp}}
\newcommand{\df}{\stackrel{\mathrm{def}}{=}}
\newtheorem{thm}{\hspace{\parindent}Theorem}[section]
\newcommand{\Pf}{{\bf Proof. }}
\newcommand{\f}{\varphi}
\newcommand{\bl}{\blacksquare}

\begin{document}

\newcommand{\vse}{\vspace{.2in}}
\numberwithin{equation}{section}

\title{On a trace formula for functions of noncommuting operators}
\author{A.B. Aleksandrov, V.V. Peller and D.S. Potapov}
\thanks{The research of the first author is partially supported by the RFBR  17-01-00607.
The publication was prepared with the support of the ``RUDN University Program 5-100''.
The research of the third author is partially supported by  ARC}

\begin{abstract}
The main result of the paper  
is that the Lifshits--Krein trace formula  cannot be generalized to the case of functions of noncommuting self-adjoint operators. To prove this, we show that for pairs $(A_1,B_1)$ and
$(A_2,B_2)$ of bounded self-adjoint operators with trace class differences $A_2-A_1$ and $B_2-B_1$, it is impossible to estimate the modulus of the trace of the difference $f(A_2,B_2)-f(A_1,B_1)$ in terms of the norm of $f$ in the Lipschitz class.
\end{abstract}

\maketitle


\

\section{\bf Introduction}
\setcounter{equation}{0}
\label{Intro}

\

The trace formula for functions of self-adjoint operators and the notion of  spectral shift function appeared in the paper \cite{L} by physicist  I.M. Lifshits in connection with problems of quantum statistics and crystals theory.
Later M.G. Krein in \cite{Kr} (see also his papers \cite{Kr2} and \cite{Kr3}) considered a significantly more general situation, extended the Lifshits result to the case of arbitrary self-adjoint operators with trace class difference, gave a mathematically rigorous definition of the spectral shift function and gave a mathematically rigorous proof of the trace formula. 

Well then, let $(A_1,A_2)$ be a pair of (not necessarily bounded) self-adjoint operators with trace class difference. We refer the reader to the books \cite{GK} and \cite{BS2} for information on the trace class $\bS_1$ and, on particular, on the libear functional trace on this class). For such a pair, there exists a unique integrable function $\bs\xi$ on the real line (called {\it the spectral shift function of the pair $(A_1,A_2)\;$}) such that 
\bay
\label{foslsso}
\trace\big(f(A_2)-f(A_1)\big)=\int_\R f'(t)\bs{\xi}(t)\,dt
\ey 
for every sufficiently nice function $f$. For example, one can take for $f$ 
a rational function with poled off the real. Note that the Krein approach is based on perturbation determin\-ants.

Next, M.G. Krein showed that the trace formula can be extended to the class of functions $f$ whose derivative is the Fourier transform of a finite Borel measure on $\R$.
He also observed that the right-hand side of the equality is well-defined for an arbitrary {\it Lipschitz function} $f$, i.e., for a function satisfying the inequality $|f(x)-f(y)|\lb\le\const|x-y|$, $x,\;y\in\R$.
In this connection he posed the question of whether trace formula \rf{foslsso} can be extended to the class of arbitrary Lipschitz function $f$. Krein also posed in \cite{Kr}  the problem to describe the maximal class of functions $f$, for which trace formula \rf{foslsso}
holds for arbitrary pairs of self-adjoint operators with trace class difference.

It turned out, however that the answer to Krein's question is negative. Indeed, Yu.B. Farforovskaya constructed in \cite{F} an example of a Lipschitz function $f$ on $\R$ and self-adjoint operators $A_1$ and $A_2$ with trace class difference such that $f(A_2)-f(A_1)\not\in\bS_1$.

As for the Krein problem of describing the maximal class of functions, for which the Lifshitz--Krein trace formula is valid, it remained open for over 60 years and was solved in \cite{Pe}. It turned out (see \cite{Pe}) that the maximal class in question coincides with the class of operator Lipschitz functions. Recall that a function $f$ on $\R$ is called {\it operator Lipschitz} if the inequality
$$
\|f(A)-f(B)\|\le\const\|A-B\|
$$
holds for arbitrary self-adjoint operators $A$ and $B$
(no matter, bounded or unbounded). We refer the reader to the survey \cite{AP} for detailed information on operator Lipschitz functions. In particular, formula \rf{foslsso} holds for functions $f$ in the homogeneous Besov class $B_{\be,1}^1(\R)$ which was proved earlier in \cite{Pe1} and \cite{Pe2}. 

Recall that it was shown in \cite{F1} that not all Lipschitz functions are operator Lipschitz, while in \cite{Mc} and \cite{K} it was established that the function $x\mapsto|x|$ is not operator Lipschitz.

Let us also mention that in \cite{BS} an attempt was undertaken to obtain trace formula \rf{foslsso} with the help of double operator integrals. Using their approach, the authors of \cite{BS} managed to prove that for each pair $(A_1,A_2)$ of self-adjoint operator with trace class difference,
there exists a finite real Borel measure $\nu$ on $\R$ (which can be called the spectral shift measure) such that
\bay
\label{pomere}
\trace\big(f(A_2)-f(A_1)\big)=\int_\R f'(t)\,d\nu(t)
\ey
for sufficiently nice functions $f$. Clearly, it follows from Krein's theorem that the measure $\nu$ is absolutely continuous with respect to Lebesgue measure $\m$ and $d\nu=\bs\xi\,d\m$. 

The program to obtain the full strength Krein theorem with the help of double operator integrals was realized in the recent paper   \cite{MNP}. To achieve this, the authors of \cite{MNP} had to consider the problem of getting a trace formula for functions of contractions. It turned out that the absolute continuity of the spectral shift measure can be deduced from the Sz.-Nagy--Foia\c s theorem on the absolute continuity of the spectral measure of the minimal unitary dilation of a completely nonunitary contraction (see \cite{SNF}).

We also mention the paper \cite{PSZ} in which a purely real method to prove the Krein theorem is given.

We would like to draw the reader’s attention to the fact that trace formula \rf{foslsso} (and even its weakened version \rf{pomere}) implies the estimate
$$
\big|\trace\big(f(A_2)-f(A_1)\big)\big|\le\const\|f\|_\Li
$$
for nice functions $f$,
where $\|f\|_\Li$ is the (semi)norm of the function $f$ in the space of Lipschitz functions. However, as follows from the result of \cite{F} mentioned above, it is impossible to replace  the modulus of the trace with the trace norm of the difference $f(A_2)-f(A_1)$ on the left-hand side of this inequality.

In this paper we work on the problem of whether one can obtain an analogue of trace formula 
\rf{foslsso} for functions of pairs of noncommuting self-adjoint operators. Such functions of operators can be defined as double operator integrals
$$
f(A,B)=\iint_{\R^2}f(s,t)\,dE_A(s)\,dE_B(t)
$$
where $E_A$ and $E_B$ are the spectral measures of the operators $A$ and $B$, see \cite{ANP}. Recall that it was shown in \cite{ANP} that under the assumption $1\le p\le2$, for functions $f$ in the homogeneous Besov class $B_{\be,1}^1(\R^2)$ the following Lipschitz type estimate
$$
\big\|f(A_2,B_2)-f(A_1,B_1)\big\|\le\const\|f\|_{B_{\be,1}^1}
\max\big\{\|A_2-A_1\|_{\bS_p},\|B_2-B_1\|_{\bS_p}\big\}
$$
holds for arbitrary pairs $(A_1,B_1)$ and $(A_2,B_2)$ of bounded self-adjoint operators such that $A_2-A_1$ and $B_2-B_1$ belong to the trace class $\bS_p$. In particular, this is true for $p=1$.

A question arizes in a natural way  of whether there is an analogue of trace formula \rf{foslsso} for functions of pairs of noncommuting operators. More precisely, let 
$(A_1,A_2)$ and $(B_1,B_2)$ be pairs of not necessarily commuting self-adjoint operators such that $A_2-A_1\in\bS_1$ and 
$B_2-B_1\in\bS_1$. Is it true that there exist integrable functions  $\bs\xi_1$ and $\bs\xi_2$ on $\R^2$ such that
\begin{align*}
\trace\big(f(A_2,B_2)-f(A_1,B_1)\big)&=
\int_{\R^2}\frac{\partial f}{\partial x}(x,y)\bs\xi_1(x,y)\,dx\,dy\\[.2cm]
&+\int_{\R^2}\frac{\partial f}{\partial y}(x,y)\bs\xi_2(x,y)\,dx\,dy
\end{align*}
for sufficiently nice functions $f$, for example, for functions 
$f$ of class 
$B_{\be,1}^1(\R^2)$?

If this is not true, is still possible to obtain a weaker result, i.e., to generalize trace formula \rf{pomere}? In other words, under the same hypotheses on the pairs of operators, do there exist finite Borel measures $\nu_1$ and $\nu_2$ such that the formula
\begin{align*}
\trace\big(f(A_2,B_2)-f(A_1,B_1)\big)=
\int_{\R^2}\frac{\partial f}{\partial x}(x,y)\,d\nu_1(x,y)
+\int_{\R^2}\frac{\partial f}{\partial y}(x,y)\,d\nu_2(x,y)
\end{align*}
holds for sufficiently nice functions $f$?

Note that each of the above formulae would imply the following estimate for the trace:
$$
\big|\trace\big(f(A_2,B_2)-f(A_1,B_1)\big)\big|\le\const\|f\|_\Li.
$$

In Section \ref{primer+} of this paper we show that all these statements for functions of noncommuting self-adjoint operators are false. 

\

\section{\bf Operator Lipschitz functions and perturbation by trace class operators}
\setcounter{equation}{0}
\label{OLyad}

\

In this section we remind properties of operator Lipschitz functions that will be used in this paper. A function $f$ defined on an interval $J$  of the real line is called {\it operator Lipschitz} if the inequality
$$
\|f(A)-f(B)\|\le\const\|A-B\|
$$
holds for arbitrary self-adjoint operators $A$ and $B$ with spectra in $J$.

First of all, we remind that operator Lipschitz functions on an on an interval $J$ 
are necessarily differentiable everywhere on $J$. This was established in \cite{JW}, see also the survey 
\cite{AP}. However, as shown in \cite{KS}, they do not have to be continuously differentiable, see also the survey \cite{AP}.

We need the following result on the behavior of functions of operators under trace class perturbations:

\medskip

{\bf Theorem on trace class perturbations.} {\it Let $J$ be an interval of the real line and let $f$ be a continuous function on $J$. The following are equivalent:

{\em(a)} $f$ is operator Lipschitz;

{\em(b)} the following inequality holds
$$
\|f(A)-f(B)\|_{\bS_1}\le\const\|A-B\|_{\bS_1}
$$
for arbitrary self-adjoint operators $A$ and $B$ with spectra in $J$ and trace class difference $A-B$.

{\em(c)} $f(A)-f(B)\in\bS_1$, whenever $A$ and $B$ are not necessarily bounded self-adjoint operators with difference $A-B$ in $\bS_1$ and spectra in $J$.}

\medskip

Moreover, if $f$ is not operator Lipschitz, then for each positive number $\e$, there exist self-adjoint operators  $A$ and $B$ such that
$A-B\in\bS_1$, $\|A-B\|_{\bS_1}<\e$ but $f(A)-f(B)\not\in\bS_1$.

We refer the reader to the survey \cite{AP}, Theorem 3.6.5, where the case $J=\R$
is considered. In the general case the proof is exactly the same.

\

\section{\bf The main result}
\setcounter{equation}{0}
\label{primer+}

\

The main purpose of this section is to show that analogues of the Lifshits--Krein trace formula for pairs of not necessarily commuting self-adjoint operators that were discussed in the introduction do not hold.

Let $j$ be a positive integer. Consider a real-valued function $\xi_j$ on $\R$ such that
$$
\xi_j(x)=\left\{\begin{array}{ll}
\big|x-2^{-j}\big|,&\big|x-2^{-j}\big|\le2^{-j-3},\\[.2cm]
0,&\big|x-2^{-j}\big|\ge2^{-j-2}
\end{array}\right.
$$
and such that $\xi_j\in C^\be(\R\setminus\{2^{-j}\})$ and $\|\xi_j\|_\Li\le\const$.
By the Johnson--Williams theorem (see \S\:\ref{OLyad}), the function $\xi_j$ is not operator Lipschitz, and so there exist self-adjoint operators $A_{1,j}$ and 
$A_{2,j}$ whose spectra are contained in $\big\{x\in\R:~\big|x-2^{-j}\big|\le2^{-j-3}\big\}$ and such that
$A_{2,j}-A_{1,j}\in\bS_1$, $\big\|A_{2,j}-A_{1,j}\big\|_{\bS_1}<2^{-j}$ but
$\xi_j\big(A_{2,j}\big)-\xi_j\big(A_{1,j}\big)\not\in\bS_1$ (see \S\:\ref{OLyad}).

It follows that we can uniformly approximate the function $\xi_j$ by an infinitely smooth real-valued function  $\psi_j$ such that
$$
\supp\psi_j\subset\big\{x\in\R:~\big|x-2^{-j}\big|\le2^{-j-3}\big\},
$$
$\|\psi_j\|_\Li\le2\|\xi_j\|_\Li$ and 
$\big\|\psi_j\big(A_{2,j}\big)-\psi_j\big(A_{1,j}\big)\big\|_{\bS_1}\ge1$.

Consider the self-adjoint operator 
$Q_j\df\psi_j\big(A_{2,j}\big)-\psi_j\big(A_{1,j}\big)$ and let $E_j$ be its spectral measure. Let $B_j$ be the self-adjoint (and unitary at the same time!) operator, which coincides with the identity operator $I$ on the subspace $E_j\big([0,\be)\big)$ and coincides with the operator $-I$ on the subspace $E_j\big((-\be,0)\big)$. Then, obviously, the operator
$
\big(\psi_j\big(A_{2,j}\big)-\psi_j\big(A_{1,j}\big)\big)B_j
$
is nonnegative and
$$
\trace\big(\big(\psi_j\big(A_{2,j}\big)-\psi_j\big(A_{1,j}\big)\big)B_j\big)
=\big\|\psi_j\big(A_{2,j}\big)-\psi_j\big(A_{1,j}\big)\big\|_{\bS_1}.
$$

Let us define now the operators $A_1$, $A_2$ and $B$ as the orthogonal sums:
$$
A_1=\bigoplus_{j\ge1}A_{1,j},\qquad
A_2=\bigoplus_{j\ge1}A_{2,j},\qquad\mbox{and}\qquad
B=\bigoplus_{j\ge1}B_j.
$$
Clearly, $A_2-A_1\in\bS_1$ and
$$
\|A_2-A_1\|_{\bS_1}=\sum_{j\ge1}\big\|A_{2,j}-A_{1,j}\big\|_{\bS_1}\le1.
$$

\begin{thm}
Let $A_1$, $A_2$ and $B$ be the self-adjoint operators defined above.
There are no complex Borel measures $\nu_1$ and $\nu_2$ such that 
the trace formula
\begin{align}
\label{fosle}
\trace\big(f(A_2,B)-f(A_1,B)\big)&=\int_{\R^2}
\left(\frac{\partial f(x,y)}{\partial x}\,d\nu_1(x,y)\right)\nonumber\\[.2cm]
&+
\int_{\R^2}\left(\frac{\partial f(x,y)}{\partial y}\,d\nu_2(x,y)\right)
\end{align}
holds for an arbitrary infinitely smooth function $f$ with compact support.
\end{thm}

\Pf For a positive integer $n$, we define the infinitely smooth function $\f_n$ on $\R$ by
$$
\f_n=\sum_{j=1}^n\psi_j,
$$
where the $\psi_j$ are the functions constructed above.
It follows easily from the definition of $\psi_j$ that $\|\f_n\|_\Li\le\const$.
Consider an infinitely smooth function $\o$ on $\R$ such that $\o(y)=y$ for $y$ in $[-1,1]$ and $\supp\o\subset[-2,2]$. Finally, we define the function $f_n$ on $\R^2$ by
$$
f_n(x,y)=\f_n(x)\o(y),\quad(x,y)\in\R^2.
$$
Clearly, the function $f_n$ is infinitely smooth, has compact support
and $\|f_n\|_{\Li}\le\const$. 

If such measures $\nu_1$ and $\nu_2$ existed, equality \rf{fosle} would imply that
$$
\big|\trace f_n(A_2,B)-f_n(A_1,B)\big|\le\const\|f_n\|_\Li\le\const.
$$
However, it is easy to see that
\begin{align*}
\trace\big(f_n(A_2,B)-f_n(A_1,B)\big)&=
\sum_{j=1}^n\trace\big(\big(\psi_j(A_{2,j}-\psi_j(A_{1,j})\big)B_j\big)\\[.2cm]
&=
\sum_{j=1}^n\big\|\psi_j\big(A_{2,j}\big)-\psi_j\big(A_{1,j}\big)\big\|_{\bS_1}
\ge n.
\end{align*}
We get a contradiction. $\bl$

\

\section{\bf Open problems}
\setcounter{equation}{0}
\label{zad}

\

{\bf1. Estimates of the moduli of traces.} We have mentioned in the introduction that in the case of functions of a single self-adjoint operator, the modulus of the trace of the difference $f(A)-f(B)$ admits a considerably stronger estimate than 
the trace norm $\|f(A)-f(B)\|_{\bS_1}$. A natural question arises of whether the same can be said in the case of functions of two noncommuting self-adjoint operators. It would be interesting to obtain a kind of an optimal estimate for
$\big|\trace\big(f(A_2,B)-f(A_1,B)\big)\big|$. This could lead to finding a version of a trace formula.

\medskip

{\bf2. A trace formula for functions of normal operators.} In \S\:\ref{primer+} of this paper we have shown that the Lifshits--Krein trace formula \rf{foslsso} cannot be generalized to the case of pairs of noncommuting self-adjoint operators. What is the situation in the case of functions ofpairs of commuting self-adjoint operators?

Clearly, it is the same as to consider functions of normal operators. 
Recall that in the paper \cite{APPS} it was shown that the inequality
$$
\big\|\big(f(N_2)-f(N_1)\big)\big\|_{\bS_p}\le\const\|f\|_{B_{\be,1}^1}\|N_2-N_1\|_{\bS_p}
$$
holds for any function $f$ of Besov class $B_{\be,1}^1(\R^2)$ and for any
 $p$ in $[1,\be]$; moreover, the constant does not depend on $p$.

The question is the following: let $N_1$ and $N_2$ be normal operators on Hilbert space. Is it true that there exist finite Borel measures  $\nu_1$ and $\nu_2$ such that the trace formula 
\begin{align*}
\trace\big(f(N_2)-f(N_1)\big)=
\int_{\R^2}\frac{\partial f}{\partial x}(x,y)\,d\nu_1(x,y)
+\int_{\R^2}\frac{\partial f}{\partial y}(x,y)\,d\nu_2(x,y)
\end{align*}
holds for sufficiently smooth functions $f$ on $\R^2$? Note that if the answer to this question is positive, then for such normal operators 
$N_1$ and $N_2$ the inequality
$$
\big|\trace\big(f(N_2)-f(N_1)\big)\big|\le\const\|f\|_\Li
$$
holds for sufficiently nice functions $f$ on $\R^2$. In any case, it would be important to find out whether under the above assumptions on can estimate $\big|\trace\big(f(N_2)-f(N_1)\big)\big|$ better than $\big\|\big(f(N_2)-f(N_1)\big)\big\|_{\bS_1}$.

If the above question has an affirmative answer, it is natural to ask the question of whether one can select measures $\nu_1$ and $\nu_2$ to be absolutely continuous with respect to Lebesgue measure on $\R^2$.

\medskip

{\bf3. A trace formula for functions of commuting self-adjoint operators.}
It was shown in \cite{NP} that if $(A_1,A_2,\cdots,A_n)$ and 
$(B_1,B_2,\cdots,B_n)$ are collections of commuting self-adjoint operators such that
$B_j-A_j\in\bS_p$, $1\le j\le n$, and $1\le p\le\be$,
then
$$
\big\|f(B_1,\cdots,B_n)-f(A_1,\cdots,A_n)\big\|_{\bS_p}
\le\const\|f\|_{B_{\be,1}^1}\max_{1\le j\le n}\|B_j-A_j\|_{\bS_p}
$$
for every function $f$ of Besov class $B_{\be,1}^1(\R^n)$.

In the case when the question posed in Subsection 2 has an affirmative answer, it would be interesting to find out whether one can generalize the Lifshits--Krein trace formula \rf{foslsso} to the case of collections of commuting self-adjoint operators.

\

\

\tiny
\noindent
\begin{tabular}{p{5cm}p{4.5cm}p{4.5cm}}
A.B. Aleksandrov & V.V. Peller& D.S. Potapov \\
St. Petersburg Department & Department of Mathematics&School of Mathematics \\
of Steklov Mathematical Institute & Michigan State University&\& Statistics  \\
of Russian Academy of Sciences & East Lansing, Michigan 48824&University of NSW\\
 27 Fontanka,  St. Petersburg&USA& Kensington, NSW 2052\\
 191023, Russia&and&Australia\\
email: alex@pdmi.ras.ru
&People’s Friendship University &email: d.potapov@unsw.edu.au\\
& of Russia (RUDN University)\\
& 6 Miklukho-Maklaya Street\\
&Moscow 117198, Russia\\
& email: peller@math.msu.edu
\end{tabular}

\end{document}